\documentclass[12pt,reqno]{amsart}

\usepackage{amsmath,amsthm,bm}
\usepackage{amssymb}
\usepackage{enumerate}
\usepackage{hyperref}

\def\ds{\displaystyle}

\def\ds{\displaystyle}

\newcommand{\R}{{\mathbb{R}}}

\newcommand{\es}{{\rm esssup\,}}

\def\XXint#1#2#3{{\setbox0=\hbox{$#1{#2#3}{\int}$ } 
\vcenter{\hbox{$#2#3$ }}\kern-.6\wd0}}

\def\B{{\mathcal B}}

\def\bA{{\bm A}}

\def\div{{\rm div\,}}

\def\bu{{\bm u}}
\def\bz{{\bm z}}

\def\ff{{\bm f}}
\def\bb{{\bm b}}
\def\diam{{\rm diam\,}}

\newtheorem{thm}{Theorem}[section]  
\newtheorem{example}[thm]{Example}

\numberwithin{equation}{section}
 
\title[Boundedness of the solutions of  parabolic systems]{Boundedness of the solutions of  a kind of nonlinear parabolic systems}
\author[E.A. Alfano]{Emilia Anna Alfano}
\address{Department of Mathematics
University of Salerno, Italy\\
 e-mail:   ealfano@unisa.it}

\author[L. Fattorusso]{Luisa Fattorusso}
\address{Department of Information Engineering,
 Infrastructure and Sustainable Energy,
Mediterranea University,
 of Reggio Calabria, Italy\\
 e-mail:  luisa.fattorusso@unirc.it}

\author[L. Softova]{Lubomira Softova}
\address{Department of Mathematics
University of Salerno, Italy\\
 e-mail:  lsoftova@unisa.it}


\begin{document}


\begin{abstract}
We deal with nonlinear systems of parabolic type satisfying componentwise structural  conditions. The nonlinear terms are Carath\'eodory maps having  controlled growth with respect to the solution and the gradient and the data are in anisotropic Lebesgue spaces. Under these assumptions we obtain essential boundedness of the weak solutions. 
\end{abstract}

\maketitle

\renewcommand{\thefootnote}{}




\section{Introduction}
This paper aims the boundedness of the weak solutions to the following nonlinear divergence form systems 
\begin{equation}\label{systemI}
{\bu}_t-\div \bA(x,t,{\bu},D{\bu})+{\bb}(x,t,{\bu},D{\bu})=0 \qquad (x,t)\in Q_T
\end{equation}
where $\Omega\subset \mathbb{R}^n,$ $n\geq 1,$ is a bounded domain, $Q_T=\Omega \times (0,T)$, $T>0,$ is a cylinder in $\R^n\times \R_+$ and the nonlinear operators
\begin{align*}
\bA(x,t,\bu,\bz)&: \  Q_T\times\R^N\times \R^{N\times n}\rightarrow\R^{N\times n} \\
\bb(x,t,\bu,\bz)&: \ Q_T\times\R^N\times \R^{N\times n}\rightarrow\R^N
\end{align*}
where $\displaystyle \bA=\{A_i^\alpha\}_{1\leq i \leq n}^{1\leq \alpha\leq N}$ and $\bb=(b^1,\cdots,b^N)$ 
are Carath\'eodory maps, that is,  these are measurable in $(x,t)\in Q_T$ for every $(\bu,\bz)\in\R^N\times  \R^{N\times n}$ and continuous in $(\bu,\bz)$ for almost all (a.a.)  $(x,t)\in Q_T.$  \\
Concerning linear divergence form equations of elliptic  type, it is well known by the result of De Giorgi and Nash \cite{DeG, Nash}  that any weak solution is locally H\"older
 continuous assuming that the leading coefficients are only $L^{\infty}.$ Almost ten years were necessary to show, via a counterexample, that this is not true for divergence form elliptic systems under the same minimal conditions on the coefficients because of the lack of maximum principle (cf. \cite{DeG2}). 
Notice that Giusti and Miranda in \cite{GM} proved that the De Giorgi-Nash result does not hold for quasi-linear systems, even if the coefficients are analytic functions.

In general, in order to obtain essential boundedness of the solutions, it is necessary to add some more restrictions on the structure of the operators considered. 
The simplest example is given by systems in diagonal form, or the so-called \textit{decoupled} systems.
\begin{example}
Let $\bu: \Omega\rightarrow \mathbb{R}^N$ be a weak solution of $\div(\bA(x,D \bu))=~0$ in $\Omega$ with 
$$
A_i^\alpha(x,D\bu)= \sum_{j=1}^n \sum_{\beta=1}^N \delta_{\alpha\, \beta}A_{ij}^{\alpha\, \beta}(x)D_ju^{\beta} 
$$
where $\delta_{\alpha\, \beta}$ is the Kronecker delta. Then each component $u^\alpha$ solves a single elliptic equation and $\sup_{\Omega}u^\alpha \leq \sup_{\partial\Omega}u^\alpha$ for each $\alpha=1,\ldots, N.$
\end{example}


As it concerns nonlinear elliptic type systems, in order to study the regularity properties of the operator, more structural conditions are required. For instance, divergence form operators
\begin{equation}\label{NonS} 
\div \bA(x,\bu,D\bu)= \bb(x,\bu,D\bu)
\end{equation}
have been studied under various componentwise conditions.

In \cite{NS} Ne\v{c}as and Star\'a analyzed quasi-linear systems which are diagonal for large values of $u^\alpha$. Precisely, for each $\alpha=1, \ldots, N$, it is assumed that
$$
0<\theta_\alpha\leq u^\alpha \Longrightarrow A_i^\alpha (x,\bu,D\bu)= \sum_{j=1}^n \sum_{\beta=1}^N  \delta_{\alpha\, \beta}A_{ij}^{\alpha\, \beta}(x,\bu)D_ju^{\beta} 
$$
where the  operator $\big\{A_{ij}^{\alpha \beta}\big\}_{1\leq i,j \leq n}^{1\leq \alpha,\beta \leq N}$ is supposed to be bounded and elliptic.  Then each solution of the system $\div\bA(x,\bu,D \bu)=0$ verifies 
$$
\sup_\Omega u^\alpha \leq \max\{\theta_\alpha ;\, \sup_{\partial\Omega} u^\alpha\}.
 $$ 

In \cite{LP} Leonetti and Petricca analyzed  nonlinear elliptic systems 
$$ 
-\div\bA(x,\bu,D\bu)+\bb(x,\bu,D\bu)=\ff(x)\quad x\in\Omega
$$
where they impose componentwise coercivity condition on the principal part and positivity of the lower order terms operator for large values of $u^\alpha$ that is there exists $\theta_\alpha>0$ such that  for each  $  u^\alpha\geq  \theta_\alpha,$ $\alpha=1,\ldots,N$ one has
\begin{equation}
\begin{cases} 
\nu |{\bz}^\alpha|^p-M_\alpha\leq \displaystyle  \sum_{i=1}^n A_i^\alpha(x,\bu,{\bz})z_i^\alpha \qquad
x\in\Omega\\[5pt] 
0\leq b^\alpha(x,\bu,{\bz}).
\end{cases}
\end{equation}
A bound for $\|\bu\|_{L^\infty(\Omega,\mathbb{R}^N)}$ has been obtained under that hypotheses when the $x$-behavior of the nonlinear terms $\bA(x,\bu,D\bu)$ and $\bb(x,\bu,D\bu)$ is controlled in Lebesgue (\cite{LP}) and Morrey (\cite{FS, PSf3, Sf1, Sf3}) spaces.

Concerning  nonlinear elliptic and parabolic  divergence form operators satisfying  quite  general  growth conditions,   there are known results on local boundedness of the solution, i.e.  in any open set $B\Subset \Omega, $  without any assumptions  on the boundary. In this context    we can mention the results of   Marcellini  et  al.  treating as the question of  local  boundedness so the Lipscitz regularity of the solutions  under $p,q$-growth  conditions  on the nonlinear operator, see     \cite{CMM1, CMM2, CMM3, Ma1, Ma2} and the references therein.

In our case we interested to obtain global regularity, up to the boundary of the cylinder $Q_T.$ For this goal, following the technique developed in \cite{LAD},  we assume  a priori that the solution is bounded on the  parabolic boundary  $\partial Q=\Omega\cup S_T.$  Further, estimating the level sets of the solution in  cylinders with small high  $Q_\tau=\Omega\times(t_1,t_2)$  with $0\leq t_1< t_2\leq T, t_2-t_1=\tau$ we can estimate the $L^\infty$ norm of $\bu$ in $Q_\tau.$   Then  covering $Q_T$ with  a finite number of such sub-cylinders and applying the global boundedness result in any of them we can  proceed up to the upper boundary $\Omega \times \{t=T\}.$  

Some more results concerning boundedness and regularity of the solutions for nonlinear elliptic and parabolic operators can be found in \cite{BPS1,BPS2,BPS3,BS1,BS2,Cm1,Cn1,CLM,FT,FT2, JS,P,Sf2}.  There is  a vast number of  regularity results, obtained via techniques of the variational calculus  that we cannot quote here and for this we direct  the reader to  the exposition works  of Mingione et  al. \cite{DMS,Mi} (see also the references therein).

As it concerns non-coercive operators, we can mention the works of Boccardo \cite{Bo2012} that treat  nonlinear elliptic equations and \cite{ADGM, AM} in case of non-coercive nonlinear equations in unbounded domain.

In what follows we use the standard  notation
\begin{itemize}
\itemsep=5pt
\item
$x=(x_1,\ldots,x_n)\in \R^n,$  $(x,t)\in \R^n\times\R_+;$\\
$\B_\rho(x)=\{y\in \R^n: |x-y|<\rho\}$ is a ball in $\R^n$ centered in $x$ and of radius $\rho>0,$ $|\B_\rho| \sim C\rho^n;$

\item
$\partial\Omega$ is the   boundary  of $\Omega,$  $\overline{\Omega}$ is the closure of $\Omega,$ and $|\Omega|$ stands for the Lebesgue measure of $\Omega;$
\item   
$Q_T=\Omega \times (0,T)$ denotes a cylinder with  base $\Omega$ and high $T>0;$\\
$S_T=\partial \Omega \times [0,T]$ is the {\it lateral boundary} of $Q_T,$ \\
  $\partial Q_T= S_T \cup \{ \Omega \times \{t=0\}\}$ stands for the  {\it parabolic  boundary} of $Q_T,$ and $\overline{Q}_T=\overline{\Omega} \times [0,T]$; 

\item $\chi_D$ is the {\it characteristic function} of any set $D\subset \R^{n+1};$

\item 
$\bu=(u^1,\ldots,u^N)$ is a vector function in $\R^N, N\geq 2$,\\
 where $u^\alpha (x,t): Q_T \rightarrow \R$, \, $\forall\, \alpha=1, \dots,N;$

\item $u^{\alpha,k}=\max\{u^\alpha(x,t)-k;0\}$ for any $k>0$ and all $(x,t)\in Q_T;$

\item
 $D_iu^\alpha= \partial u^\alpha /  \partial x_i $ and $D\bu=\{D_iu^\alpha\}_{i\leq n}^{\alpha\leq N}  \in \mathbb{R}^{N\times n}$  is the vector of the space derivatives of $\bu,$ $D_tu^\alpha=u_t^\alpha = \partial u^\alpha/\partial t;$ 

 
\item
 for any vector-function $\ff\in L^p(Q,\R^N), p>1$ we write $\|\ff\|_{p;Q}$ instead of $\|\ff\|_{L^p(Q;\R^N)}.$
\end{itemize}
The  letter  $C$ denotes  a positive constant depending on known quantities  that can vary  from one occurrence to another.  In addition  the summation convention  on the repeated  indexes is adopted.

\section{Definitions and auxiliary results}
\label{sec2a}
We say that the domain  $\Omega$  satisfies the  (A) condition if  there exists a  positive constant $A_\Omega<1$ such that, for each ball $\B_\rho$ centered in some point   $x\in \partial \Omega$   and $\rho\leq \diam \Omega$, there  holds
\begin{equation}\tag{A}\label{Acond}
|\Omega\cap \B_\rho|\leq (1-A_\Omega) |\B_\rho|\,.
\end{equation} 
For the  function spaces that we are going to use  we adopt  the notations introduced in \cite{LAD}.
\begin{enumerate}
\item[1.] We say that  ${\bu}=(u^1, \ldots, u^N):Q_T \rightarrow \mathbb{R}^N$ belongs to $L^2(Q_T; \R^N)$ if $u^{\alpha} \in L^2(Q_T)$ for each $\alpha=1,\ldots,N$ and 
$$
\| {\bu}\|_{2; Q_T}= \left( \int_{Q_T} |{\bu}|^2\, dxdt\right)^{1\over 2}= \left( \int_{0}^T \int_{\Omega} \sum_{\alpha=1}^N |u^{\alpha} (x,t)|^2 \, dxdt \right)^{1\over 2}.
$$

\item[2.] For any $q,r \geq 1$ we consider the  anisotropic Lebesgue space $L^{q,r}(Q;\mathbb{R}^N)$ equipped with mixed norm 
$$
\|{\bu}\|_{q,r;Q_T}=\left[ \int_0^T \left( \int_{\Omega} \sum_{\alpha=1}^N |u^{\alpha} (x,t)|^q \, dx\right)^{r\over q} dt \right]^{1\over r}.
$$
In the particular case $q=2$, $r=\infty$, we have 
$$
\|{\bu}\|_{2,\infty;Q_T}=\underset{t\in[0,T]}\es \|{\bu}(\, \cdot\, ,t)\|_{2;\Omega} 
$$
and for  $q=r$ we obtain the  parabolic Lebesgue space  $ L^q(Q;\Bbb{R}^N).$ 

\item[3.]  $W_2^{1,0}(Q_T; \R^N)$  and   $W_2^{1,1}(Q_T; \R^N)$  are the Banach spaces  consisting     of all  ${\bu} \in L^2(Q_T;\R^N)$ for which the following norms are finite 
\begin{align*}
\|{\bu}\|_{W_2^{1,0}(Q_T)}=\, & \|{\bu}\|_{2;Q_T}+\|D {\bu}\|_{2;Q_T}\\ 
\|{\bu}\|_{W_2^{1,1}(Q_T )}=\, & \|{\bu}\|_{2;Q_T}+\|D {\bu}\|_{2;Q_T} + \|{\bu}_t\|_{2;Q_T}.
\end{align*}
 The space $\overset{\circ}{W}{}^{1,1}_2(Q_T;\R^N)$ is the closure of $C_0^{\infty}(Q_T; \mathbb{R}^N)$ with respect to the last norm. 

\item[4.] $V_2(Q_T; \mathbb{R}^N)$ is a subspace of $W_2^{1,0}(Q_T;\mathbb{R}^N)$ such that
\begin{equation}\label{norm5}  
\|{\bu}\|_{V_2(Q_T)}= \|{\bu}\|_{2,\infty;Q_T} + \|D{\bu}\|_{2;Q_T}<\infty. 
\end{equation}
Namely,
$$
 V_2(Q_T; \R^N)=L^{\infty}(0,T;L^2(\Omega;\R^N) ) \cap  L^2(0,T;W_2^1(\Omega;\R^N) )
$$ 
endowed with the norm \eqref{norm5}. 
\item[5.]
$V_2^{1,0}(Q_T; \R^N)$ consists of all ${\bu} \in V_2(Q_T; \mathbb{R}^N)$ continuous in $t$ with respect to the norm of $L^2(Q_T;\R^N)$, that is 
\begin{equation}\label{norm6}   
\lim_{\Delta t \rightarrow 0} \|{\bu}(\cdot, t+\Delta t)- {\bu}(\cdot, t)\|_{2;\Omega}=0
\end{equation}
endowed with the norm \eqref{norm5}.
As usual, we denote by 
 $\overset{\circ}{V}{}_2(Q_T;\R^N)$ the closure of $C_0^\infty(Q_T;\R^n)$ with respect to the norm \eqref{norm5}.
\end{enumerate}

For reader's convenience, we recall the \textit{H\"older inequality  on anisotropic Lebesgue spaces} and some of its modifications  (cf. \cite[Ch. II]{LAD}). Precisely, for any exponents $q\geq q_1\geq 1,$  $r\geq r_1\geq 1,$ there hold
\begin{align}\label{Holder} 
\Big| \int_{Q_T} u(x,t)\, v(x,t)\, dxdt \Big| & \leq \|{u}\|_{q,r;Q_T}\|{v}\|_{\frac{q}{q-1}, \frac{r}{r-1};Q_T}\\[5pt]
\label{1.7LSU}
\|uv\|_{q_1,r_1;Q_T} &  \leq \|u\|_{q,r;Q_T}  \|v\|_{\frac{qq_1}{q-q_1},\frac{rr_1}{r-r_1};Q_T}.
\end{align}

Now, we give also  some properties of the spaces $  V_2(Q_T)$ and $ \overset{\circ}{V}{}_2(Q_T).$ Let  $u \in \overset{\circ}{V}{}_2(Q_T)$, then   from  the Gagliardo--Nirenberg interpolation inequality for the   $\overset{\circ}{W}{}^1_2(\Omega)$-functions there exists a constant $C_q=C_q(n,q)$ such that
$$
\|u(\cdot,t)\|_{q,\Omega}\leq C_q \|u(\cdot,t)\|^{1-\alpha}_{2,\Omega} \|Du(\cdot,t)\|^\alpha_{2,\Omega}\quad \text{ for a.a. } t\in[0,T]
$$  
with $\alpha=\frac{n}2-\frac{n}q,$ and where $q$ is as in \eqref{intq,r}. Hence, we have the following imbedding type inequality in  $\overset{\circ}{V}{}_2(Q_T)$
$$
 \|u\|_{q,r;Q_T} \leq C_q  \|u\|_{2,\infty;Q_T}^{1- \alpha} \, \|Du\|_{2,\alpha r;Q_T}^\alpha  . 
 $$
Taking  $\alpha={2/ r},$ where $r$ is as in \eqref{intq,r}, we get  a series of inequalities and the first one is 
\begin{equation}\label{3.2}
\|u\|_{q,r;Q_T}\leq C_q \|u\|_{2, \infty;Q_T}^{1-{2\over r}} \, \|Du\|_{2;Q_T}^{2 \over r}
\end{equation}
where $\ds {1\over r}+ {n\over 2q}={n \over 4},$ and $r$ and $ q $ satisfy 
\begin{equation}\label{intq,r}
\begin{cases} 
r \in [2,\infty] \quad &q \in [2, {2n \over n-2}], \quad  \mbox{ for }n\geq 3 \\[4pt]
r \in (2,\infty) \quad &q \in [2, \infty),\quad \ \, \mbox{ for }n=2 \\[4pt]
r \in [4,\infty] \quad &q \in [2, \infty],\quad \ \  \mbox{ for }n=1.
\end{cases}
\end{equation}

Estimating the right-hand  side of \eqref{3.2} by Young's inequality
$$
 ab\leq \frac2r\ a^{r\over 2}+\frac{r-2}r\ b^{r \over r-2}, \qquad r>2 
$$
we obtain the second inequality
\begin{equation}\label{a1}
 \|u\|_{q,r;Q_T}\leq C_q {\frac2r} \, \|Du\|_{2;Q_T}+ C_q \frac{r-2}r \, \|u\|_{2, \infty;Q_T} \, \leq \beta \, \|u\|_{V_2(Q_T)}
\end{equation}
where $\beta=\beta(n,q,r)$ and $q$ and $r$ are as in \eqref{intq,r}.

One more notion that we need is that of {\it the Steklov average.}
For a function $\zeta(x,t) \in \overset{\circ}{W}{}^{1,1}_2(Q_{[-h,T]})$ that is zero for $t\leq 0$ and $t\geq T-h,$  we define \textit{the Steklov average in the future} as 
\begin{equation}\label{future}
\zeta_h (x,t)={1\over h}\int_t^{t+h} \zeta(x,\tau)\, d\tau 
\end{equation}
and {\it the Steklov average in the past} as 
\begin{equation}\label{past}
\zeta_{\overline{h}}(x,t)={1\over h}\int_{t-h}^t \zeta(x,\tau)\, d\tau. 
\end{equation}

We notice that, for each fixed $h$, $\zeta_h (x,t)$ is defined in $\overline{Q}_{T -h} =\overline{\Omega}\times [0,T-h]$ and $\zeta_{\overline{h}}(x,t) $ is defined in $\overline{Q}_{[-h,T]} =\overline{\Omega}\times [-h,T],$
 where $0<h<T.$ 
It is well known that if $\zeta \in L^{q,r}(Q_T)$ with $q,r\geq 1$, then $\zeta_h$ approximates $\zeta$ with respect to the norm  in $L^{q,r}.$ 
Furthermore, if $\zeta \in \overset{\circ}{V}{}_2(Q_{T})$, then $\zeta_h \in \overset{\circ}{W}{}^{1,1}_2(Q_{T-h})$ and  $\zeta_h$ approximates $\zeta$ with respect to the norm in  $V_2$ (cf. \cite[Ch. II]{LAD}). 
Moreover, if we take  $\zeta \in \overset{\circ}{V}{}_2(Q_T),$ then $\zeta_{\overline{h}} \in \overset{\circ}{W}{}^{1,1}_2(Q_{[h,T]}).$

We recall that the Steklov averages satisfy these two properties: 
\begin{equation}\label{rel}
\begin{split}
\left(\zeta_{\overline{h}}\right)_{t} &=  \left(\zeta_t\right)_{\overline{h}}\\
\int_0^T u(x,t)\, \zeta_{\overline{h}} (x,t)\, dt&=\int_0^{T-h} u_h(x,t)\, \zeta(x,t)\, dt.
\end{split}
\end{equation}

The last equation holds for any square summable function, such that $\zeta(x,t)=0$ for $t \leq 0$ and for  $t \geq T-h.$

The following result permits to give a total estimate for the maximum of the modulus of the solutions in the whole domain of the definition. 
Let $u$ in $V_2(Q_T)$. Denote by 
$u^k=\max\{u(x,t)-k;0\},$ $k>0 $ for all $(x,t)\in Q_T$ and define the set 
$$
A^k(t)=\{x\in\Omega:\, u(x,t)>k, \, t \in[0,T]\}.
$$
\begin{thm}[\cite{LAD}, Ch. II, Theorem 6.1]\label{T}
Suppose that $\|u\|_{\infty;S_T}\leq M_0$, $M_0 \geq 0$,   and 
\begin{equation}\label{t1}
\|u^{k}\|_{V_{2}^{1,0}(Q_{T})}\leq C k \mu^{{1+\varkappa}\over r}(k) 
\end{equation}
hold for $k\geq M_0$, where $\mu(k)=\int_0^T |A^k(t)|^{r\over q} \,dt$, $q$ and $r$ as in \eqref{intq,r} and $0<\varkappa<1$. Then 
$$
\|u\|_{\infty;Q_{T}}\leq 2M_0\left[1+2^{{2\over \varkappa}+{1\over \varkappa^2}} (\beta c)^{1+{1\over \varkappa}} T^{{1+\varkappa \over r}}
 |\Omega|^{{1+\varkappa \over q}}\right]
$$
where $\beta$ is the constant from \eqref{a1} and $C$ is an arbitrary constant. 

\end{thm}

\section{Statement of the problem}\label{sec2}

Recall that we are going to study the boundedness of the  weak  solutions to the system 
\begin{equation}\label{system}
u^\alpha_t-\sum_{i=1}^n D_i\big(A_i^\alpha(x,t,{\bu},D{\bu})\big)+{b^\alpha}(x,t,{\bu},D{\bu})=0 
\end{equation}
in $Q_T$ with $\alpha=1,\cdots,N$.
In order to obtain our results we impose the following   assumptions regarding the structure and the asymptotic behaviour  of the nonlinear operators  $\bA$ and $\bb.$ 

\begin{itemize}
\itemsep=3pt
\item[$(H_1)$]  {\it Controlled growth conditions:} there exist a positive constant $\Lambda$  and  functions  $\varphi_1\in L^2(Q_T)$ and $\varphi_2 \in L^{q_0, r_0}(Q_T),$ with  $r_0,q_0\geq 1,$ ${1\over r_0}+ {n\over 2q_0}=1+{n \over 4},$ such that 
\begin{align}\label{1.4}
|{\bA}(x,t,\bu,\bz)|&\leq\Lambda\left(\varphi_1(x,t)+|{\bu}|^{\frac{n+2}{n}}+|{\bz}|\right)\\
\label{1.5}
|{\bb}(x,t,\bu,\bz)|&\leq\Lambda\left(\varphi_2(x,t)+|\bu|^{1+\frac{4}{n}}+|\bz|^{\frac{n+4}{n+2}}\right)
\end{align}
for a.a. $(x,t)\in Q_T$ and all $(\bu,\bz)\in \R^N\times\R^{N\times n}.$

\item[$(H_2)$]   {\it  Componentwise coercivity  of the differential operator:}    there exist positive constants $\nu, \mu$  and $k_0$ such that for all $\alpha\in\{1,\ldots,N  \}$ and for $|{u^\alpha}|\geq k_0$, it holds
\begin{equation} \label{cond1}
\sum_{i=1}^n A^\alpha_i(x,t,\bu,\bz)z^\alpha_i\geq \nu |{\bz}^\alpha|^2-\mu|u^\alpha|^\delta- {|u^\alpha|}^2 \psi(x,t)
\end{equation}
for a.a. $(x,t)\in Q_T.$    The function $\psi\in L^{q_2,r_2}(Q_T)$ is as in \eqref{1.8}.

\item[$(H_3)$]   {\it  Componentwise  sign condition:} for $|{u^{\alpha}}|\geq k_0$  we have
\begin{equation} \label{cond2}
-b^\alpha(x,t,{\bu},{\bz})u^\alpha\leq |u^\alpha|^2 \psi(x,t) +\mu |u^\alpha|^{\delta}+\Lambda |\bz^\alpha|^2 
\end{equation}
for a.a. $(x,t)\in Q_T$ and for all ${\bz}\in \R^{N\times n}.$
\end{itemize}

The vector function  $\bu\in V^{1,0}_2(Q_T;\R^N) $  is a {\it weak solution} of \eqref{system} if, for any  $0\leq t_0 \leq t_1\leq T$ and each $\alpha=1,\ldots,N$ it holds 
\begin{equation}\label{weak}
\begin{split}
-\int_{t_0}^{t_1}  \int_\Omega & u^\alpha(x,t) \eta_t^\alpha(x,t) \,  dx dt \\ 
&+ \int_{t_0}^{t_1}\int_\Omega \Big\{ \sum_{i=1}^n A^\alpha_i(x,t,\bu,D \bu)D_i\eta^\alpha(x,t) \\
&+b^\alpha(x,t,{\bu},D{\bu})\eta^\alpha(x,t)\Big\}\,dxdt\\
&+\int_\Omega u^\alpha(x,t)\eta^\alpha(x,t)\Big\vert_{t_0}^{t_1}\,dx =0 
\end{split}
\end{equation}
for any $ \bm{ \eta}(x,t)\in {\overset{\circ}{W}}{}^{1,1}_2(Q_T;\R^N).$
The definition  \eqref{weak} has sense if the integrals are convergent and this is possible because of the {\it growth conditions} $(H_1).$  
Now, we take as  test function in \eqref{weak}  the Steklov average 
$$
\zeta^\alpha_{\overline{h}}= {1\over h} \int_{t-h}^t \zeta^\alpha(x,\tau)\, d\tau
$$
where $\zeta^\alpha \in \overset{\circ}{V}{}_2^{1,0}(Q_{T-h}).$ 
Then,  \eqref{weak} becomes 
\begin{equation}\label{weak1}
\begin{split}
-\int_{t_0}^{t_1} \int_\Omega & u^\alpha(x,t) ({\zeta_{\overline{h}}^\alpha})_{t}(x,t)\, dx dt \\ 
&+  \int_{t_0}^{t_1} \int_\Omega \Big\{ \sum_{i=1}^n A_i^\alpha(x,t,\bu,D \bu)\, D_i\, \zeta_{\overline{h}}^\alpha(x,t) \\
&+ b^\alpha(x,t,{\bu},D{\bu})\, \zeta_{\overline{h}}^\alpha(x,t)\Big\}\,dxdt \\ 
& + \int_\Omega u^\alpha(x,t)\zeta^\alpha_{\overline{h}}(x,t)\bigg\vert_{t_0}^{t_1}\, dx =0
\end{split}
\end{equation}
where $0\leq t_0\leq t_1\leq T-h.$
By integration by parts and \eqref{rel}  the first term in \eqref{weak1} becomes

\begin{align*}
-\int_{t_0}^{t_1}\!\! \int_\Omega& u^\alpha(x,t) \left(\zeta_{\overline{h}}^\alpha\right)_t
(x,t) \, dx dt = -\int_{t_0}^{t_1}\!\! \int_\Omega u^\alpha(x,t) \left(\zeta_t^\alpha\right)_{\overline{h}}
(x,t) \, dx dt\\
&=-\int_{t_0}^{t_1}\!\! \int_\Omega u^\alpha_h(x,t) \zeta_t^\alpha (x,t)\, dx dt\\ 
&=- \int_\Omega u^\alpha_h(x,t) \zeta^\alpha(x,t)\bigg\vert_{t_0}^{t_1}\,dx +\int_{t_0}^{t_1}\!\! \int_\Omega {u_{ht}^\alpha}(x,t)\zeta^\alpha(x,t) \, dx dt \\ 
&=- \int_\Omega u^\alpha(x,t) \zeta_{\overline{h}}^\alpha(x,t)\bigg\vert_{t_0}^{t_1} \, dx +\int_{t_0}^{t_1}\!\! \int_\Omega {u_{ht}^\alpha}(x,t)\zeta^\alpha(x,t) \, dx dt
\end{align*}
where $u_h^\alpha$ is the Steklov averege of $u^\alpha$  (see  \eqref{future}). Hence,  
\begin{equation}\label{weak2} 
\begin{split}
\int_{t_0}^{t_1}\!\! \int_\Omega & \Big[ u_{ht}^\alpha (x,t) \zeta^\alpha (x,t) + \sum_{i=1}^n A_i^\alpha (x,t,\bu,D \bu) D_i \zeta_{\overline{h}}^\alpha(x,t) \\
& + b^\alpha(x,t,\bu,D\bu) \zeta_{\overline{h}}^\alpha(x,t)\Big]\, dxdt=0. 
\end{split}
\end{equation} 
For some positive constant $k$ we take
\begin{equation}\label{eq-k}
 \zeta^\alpha(x,t)= \xi^2(x,t)\max\{u_h^\alpha(x,t)-k;0\}=: \xi^2 u_h^{\alpha,k}
\end{equation}
 where $\xi(x,t)$ is an arbitrary nonnegative continuous piecewise-smooth function that is equal to zero on the lateral boundary $S_T.$  Observing  that 
$$
D_t u^{\alpha,k}=D_t u^\alpha, \quad D_i u^{\alpha,k}=D_i u^\alpha
$$
we transform the first term in \eqref{weak2} as follows 
\begin{align*}
\int_{t_0}^{t_1} \int_\Omega u_{h t}^\alpha\, u_h^{\alpha,k}\, \xi^2\, dxdt&= {1\over 2}\int_\Omega \big( u_h^{\alpha,k}\big)^2\xi^2(x,t)
\Big\vert_{t_0}^{t_1} \, dx\\
&- \int_{t_0}^{t_1} \int_\Omega\big(u_h^{\alpha,k}\big)^2 \xi(x,t) \xi_t(x,t)\, dxdt.  
\end{align*}
Now, we pass to the limit in \eqref{weak2} as $h \rightarrow 0$, obtaining 
\begin{equation}\label{1.9LAD}
\begin{split}
{1\over 2}\, &\|u^{\alpha,k}(\, \cdot , t)\xi(\, \cdot,t)\|_{2, \Omega}^2\Big\vert_{t_0}^{t_1}  + \int_{t_0}^{t_1} \int_\Omega\Big[- \left(u^{\alpha,k}\right)^2 \xi(x,t)\, \xi_t(x,t) \\
&+ \sum_{i=1}^n A_i^\alpha (x,t,\bu,D \bu)\, D_i(u^{\alpha,k}\, \xi^2(x,t))\\
 &+ b^\alpha(x,t,{\bu},D{\bu})\, u^{\alpha,k}\, \xi^2(x,t) \Big]\, dxdt=0.
\end{split}
\end{equation}
The validity of the limit in all of the terms follows  by the properties of the Steklov average (see \cite{LAD}). 

In what follows, we assume that the weak solution  $\bu$ of \eqref{system}
belongs to $ L^{q_1,r_1}(Q_T,\R^N)$ where the exponents $q_1,r_1 >1$  and the other parameters to be used in the sequel, satisfy
\begin{equation} \label{1.8}
\begin{cases}
\text{ for } n \geq 2 & \\[5pt]
(\delta-2)\left(\dfrac{1}{r_1}+\dfrac{n}{2q_1}\right)<1, \  & \delta>2,\quad \dfrac{1}{r_2}+\dfrac{n}{2q_2}<1,\\[10pt]
 \dfrac{q_1}{\delta-2},\,  q_2\in \left(\dfrac{n}{2},\infty\right],  &\dfrac{r_1}{\delta-2}, \, r_2\in \left(1,\infty\right]\\[10pt]
  \psi\in L^{q_2,r_2}(Q_T)  & q_2\in\left(\dfrac{n}{2},\infty\right],\, r_2\in(1,\infty]\\[10pt]
  \text{ for } n=1&\\[5pt]
\dfrac{q_1}{\delta-2}, \,  q_2\in [1,\infty],  &\dfrac{r_1}{\delta-2}, \, r_2\in (1,\infty)\\[10pt]
  \psi\in L^{q_2,r_2}(Q_T) &  q_2\in[1,\infty],\, r_2\in(1,\infty).
\end{cases}
\end{equation}

Suppose that the solution $\bu$ is essentially bounded on the parabolic boundary $\partial Q_T$ with some constant $M_0>0,$ that is  $\|\bu\|_{\infty;\, \partial Q_T}\leq M_0;$ we take $k \geq \max\{k_0, M_0\},$ where $k_0>0$ is the constant from $(H_2),$  and observe that $u^{\alpha,k}$ is equal to zero on the parabolic  boundary $\partial Q_T$. Thus, in \eqref{1.9LAD} we can choose $t_0=0$ and $\xi =1,$ and hence   $\xi_t=0,$ that gives 
\begin{equation}
\label{2.4}
\begin{split}
\frac{1}{2}\|{ u}^{\alpha,k}(\, \cdot, t_1)\|^{2}_{2,A^{\alpha,k}(t_1)}& +\int_{Q^{\alpha,k}_{t_1}}\Big[ \sum_{i=1}^n A^\alpha_i(x,t,{\bu},D{\bu})\, D_i u^{\alpha,k} \\
&  +b^\alpha(x,t,{\bu},D{\bu})\, u^{\alpha,k}\Big]\,dx dt =0
\end{split}
\end{equation}
where the integration is taken on the sets 
\begin{align*}
A^{\alpha,k}(t_1)&=\{x\in\Omega: u^\alpha(x,t_1)>k,\; t_1 \in [0,T]\} \\
Q^{\alpha,k}_{t_1} &=\{(x,t)\in\Omega\times [0,t_1]: u^\alpha(x,t)>k\}. 
\end{align*}

\section{Main result}
The conditions imposed above permit us to obtain a componentwise maximum principle and, hence, also essential boundedness of the solution of \eqref{system}.
\begin{thm}
\label{T2.1}
Let $\bu \in V_2^{1,0}(Q_T, \R^{N})\cap L^{q_1,r_1}(Q_T, \R^{N})$ be a weak solution  of the system \eqref{system} which is essentially bounded on the boundary $\partial Q_T.$
Assume that the terms
 ${\bA}(x,t,{\bu},{\bz})$ and ${\bb}(x,t,{\bu},{\bz})$ satisfy conditions \eqref{1.4}-\eqref{cond2} with $\nu - \Lambda >0.$ Then 
\begin{equation}
\label{2.1}
\|u^\alpha\|_{\infty;Q_{T}}\leq  K 
\end{equation}
with a constant $K$ depending on known quantities and   $\|{ u^\alpha}\|_{q_1,r_1;Q_T}.$  
If  we take  $\mu=0$  in  \eqref{cond1}-\eqref{cond2},  then  the constant $K$ does not depend on the last norm. 
\end{thm}

\proof
By   \eqref{2.4}, taking into account $\displaystyle {u^{\alpha,k} \over u^\alpha}  <1$,  \eqref{cond1} and \eqref{cond2}, we get
\begin{equation}
\label{2.5}
\begin{split}
\frac{1}{2} & \|u^{\alpha,k}\|^{2}_{2,A^{\alpha,k}(t_1)}+\int_{Q_{t_1}^{\alpha,k}}\left(\nu|Du^{\alpha}|^2-\mu|{u^\alpha}|^{\delta}-{ |u^\alpha|}^2\psi\right)\,dx dt \\
& \leq-\int_{Q_{t_1}^{\alpha,k}}b^{\alpha}(x,t,{\bu},D{\bu})\, u^{\alpha}(x,t)\displaystyle\frac{u^{\alpha,k}(x,t)}{u^{\alpha}(x,t)}\,dx dt \\
& \leq\int_{Q_{t_1}^{\alpha,k}}\left(|u^{\alpha}|^2\psi +\mu|u^{\alpha}|^{\delta}+  \Lambda|Du^{\alpha}|^2\right)\,dxdt. 
\end{split}
\end{equation}
Hence 
\begin{equation}\label{2.5*} 
\begin{split}
\frac{1}{2}&\|u^{\alpha,k}\|^{2}_{2,A^{\alpha,k}(t_1)}+(\nu-\Lambda)\int_{Q_{t_1}^{\alpha,k}} |Du^{\alpha}|^2\,dxdt \\
& \leq 2 \int_{Q_{t_1}^{\alpha,k}} \left( \mu |u^\alpha|^\delta+ |u^{\alpha}|^2\psi \right)\,dxdt.
\end{split}
\end{equation}
Keeping in mind that $\nu-\Lambda>0$, we obtain 
\begin{equation}\label{2.5**}
\|u^{\alpha,k}\|^2_{2;A^{\alpha,k}(t_1)} + \|Du^\alpha\|^2_{ 2;Q^{\alpha,k}_{t_1}}\leq C \int_{Q^{\alpha,k}_{t_1}} \left( |u^\alpha|^\delta + |u^\alpha|^2 \psi\right)\, dxdt
\end{equation}
where $C$ depends on $\mu$ and $\nu - \Lambda.$
Let us estimate the integrals in the right-hand side. For the first one, using the inequalities in \S \ref{sec2a}, we obtain 

\begin{equation}\label{2.7}
\begin{split}
I_1&=\int_{Q_{t_1}^{\alpha,k}}|u^{\alpha}|^{\delta-2}\cdot|u^\alpha |^2\,dxdt\\
& \leq\left[\int_0^{t_1}\left(\int_{A^{\alpha,k}(t_1)}|u^\alpha|^{q_1}dx\right)^{\frac{r_1}{q_1}}dt\right]^{\frac{\delta-2}{r_1}} \\
& \quad\times \left[
\int_0^{t_1}\left(\int_{A^{\alpha,k}(t_1)}|u^{\alpha}|^{\overline{q}_1}dx\right)^{\frac{\overline{r}_1}{\overline{q}_1}}dt
\right]^{\frac{2}{\overline{r}_1}}  \\
& =\|u^\alpha\|^{\delta-2}_{q_1,r_1,Q_{t_1}^{\alpha,k}}\|u^\alpha\|^{2}_{\overline{q}_1,\overline{r}_1,Q_{t_1}^{\alpha,k}} \\
&\leq 2\|u^\alpha\|^{\delta-2}_{q_1,r_1,Q_{t_1}^{\alpha,k}}\biggl[\| 
u^\alpha -k\|^{2}_{\overline{q}_1,\overline{r}_1,Q_{t_1}^{\alpha,k}}\\
  &\qquad \qquad+k^2\left(\int_0^{t_1}( |A^{\alpha,k}(t)|)^{\frac{\overline{r}_1}{\overline{q}_1}}\,dt\right)^{\frac{2}{\overline{r}_1}}\bigg]
\end{split}
\end{equation}
where $\displaystyle \overline{q}_1={2q_1 \over q_1-(\delta-2)}$, $\displaystyle \overline{r}_1={2r_1 \over r_1-(\delta-2)}$. We notice that 
\begin{equation}\label{int}
\begin{cases} 
\displaystyle
2\leq \overline{r}_1< \infty, \quad 2\leq \overline{q}_1<{2n \over n-2} \quad  &\mbox{ for }n\geq 2 \\
\displaystyle
2< \overline{r}_1< \infty, \quad 2\leq \overline{q}_1\leq \infty \quad  &\mbox{ for }n=1 
\end{cases}
\end{equation}
and 
$\displaystyle {1 \over \overline{r}_1} + {n \over 2\overline{q}_1 }>{n \over 4}.$ If $n>2$ we can find an opportune $0<\varkappa_1<1,$ such that $\displaystyle {1 \over \overline{r}_1} + {n \over 2\overline{q}_1 }={n \over 4}(1+\varkappa_1).$ In case of $n=2 $ we consider $\overline{q}_1\in [2,\infty).$
Introduce
$$
\widehat{q}_1= \overline{q}_1 \left(1+\varkappa_1\right) \quad \text{and} \quad
\widehat{r}_1= \overline{r}_1 \left(1+\varkappa_1\right).
$$
Direct calculations show that $\widehat{q}_1$ and $\widehat{r}_1$ verify  \eqref{intq,r}
and by  \eqref{1.7LSU} we obtain
\begin{equation}\label{n}
 \begin{split}
 \|u^{\alpha,k}\|^2_{\overline{q}_1, \overline{r}_1;Q_{t_1}^{\alpha,k}}&\leq \|u^{\alpha,k}\|^2_ {\widehat{q}_1, \widehat{r}_1;Q_{t_1}^{\alpha,k}} \, \|\chi_{Q_{t_1}^{\alpha,k}}\|^2_{{\widehat{q}_1\over \varkappa_1},{\widehat{r}_1\over \varkappa_1};Q^{\alpha,k}_{t_1}}\\ 
& \leq \|u^{\alpha,k}\|^2_ {\widehat{q}_1, \widehat{r}_1;Q_{t_1}^{\alpha,k}}\, \left(t_1^{{1\over \widehat{r}_1}}  |\Omega|^{{1\over\widehat{q}_1}} \right)^{2\varkappa_1}.
 \end{split}
\end{equation}

In the case $q_1=\delta - 2,$ which is possible if $n=1,$ we obtain $\overline{q}_1=\infty$ and hence also $\widehat{q}_1=\infty.$ Then \eqref{n} becames
\begin{equation}\label{n1}
 \|u^{\alpha,k}\|^2_{\infty, \overline{r}_1;Q_{t_1}^{\alpha,k}}\leq  \|u^{\alpha,k}\|^2_ {\infty, \widehat{r}_1;Q_{t_1}^{\alpha,k}}\, t_1^{{2\varkappa_1\over \widehat{r}_1}}.
\end{equation}

This permits us  to apply  \eqref{a1}. In the case $n\geq 2$ we obtain 
\begin{equation}\label{2.7*}
\begin{split}
I_1& \leq 2\|u^\alpha\|^{\delta-2}_{q_1,r_1,Q_{t_1}^{\alpha,k}}\Biggl[\beta^2\|u^{\alpha,k}\|_{V_2^{1,0}(Q_{t_1}^{\alpha,k})}^2 \\
& \quad\times\left(\int_0^{t_1}(|A^{\alpha,k}(t)|)^{\frac{\widehat{r}_1}{\widehat{q}_1}}dt\right)^{\frac{2}{\overline{r}_1}-\frac{2}{\widehat{r}_1}} \\
& \left.+k^2\left(\int_0^{t_1}(|A^{\alpha,k}(t)|)^{\frac{\widehat{r}_1}{\widehat{q}_1}}dt\right)^{\frac{2(1+\varkappa_1)}{\widehat{r}_1}}\right] \\
& \leq 2 \|u^\alpha\|^{\delta-2}_{q_1,r_1,Q_{t_1}^{\alpha,k}} \Bigg\{\beta^2 \|u^{\alpha,k}\|_{V_2^{1,0}(Q_{t_1}^{\alpha,k})}^2   \big(t_1^{\frac{1}{\widehat{r}_1}} |\Omega |^{\frac{1}{\widehat{q}_1}} \big)^{2\varkappa_1}\\ 
& + k^2\left(\int_0^{t_1}\,|A^{\alpha,k}(t)|^{\frac{\widehat{r}_1}{\widehat{q}_1}}\,dt\right)^{\frac{2(\varkappa_1+1)}{\widehat{r}_1}}\Bigg\}
\end{split}
\end{equation}
where 
$\beta=\beta (n, \overline{q}_1, \overline{r}_1)$
 and $\displaystyle {\widehat{q}_1 \over \overline{q}_1}={\widehat{r}_1 \over \overline{r}_1}=\varkappa_1 +1.$

If $q_1=\delta-2,$ in view of \eqref{n1}, we have 
\begin{equation}
\begin{split}
I_1&\leq 2 \|u^\alpha\|^{\delta-2}_{\delta-2,r_1,Q_{t_1}^{\alpha,k}} \left(\beta^2 \|u^{\alpha,k}\|_{V_2^{1,0}(Q_{t_1}^{\alpha,k})}^2   t_1^{\frac{2\varkappa_1}{\widehat{r}_1}} + k^2\, t_1^{\frac{2(\varkappa_1+1)}{\widehat{r}_1}}\right).  
\end{split}
\end{equation}
We trate $I_2$ analogously. By the H\"older inequality \eqref{Holder}, we obtain
\begin{equation}
\label{2.8}
\begin{split}
I_2=\int_{Q_{t_1}^{\alpha,k}} |u^\alpha|^2 \psi \, dxdt
& \leq\|\psi\|_{q_2,r_2,Q_{t_1}^{\alpha,k}}\, \|u^{\alpha}\|_{\overline{q}_2,\overline{r}_2,Q^{\alpha,k}_{t_1}}^2
\end{split}
\end{equation}
where $\ds \overline{q}_2=\frac{2q_2}{q_2-1}$, $\ds \overline{r}_2=\frac{2r_2}{r_2-1}$. We notice that 
\begin{equation*}
\begin{cases} 
\displaystyle {1 \over \overline{r}_2} + {n \over 2\overline{q}_2 }>{n \over 4} &\\[6pt]
2\leq \overline{q}_2<\dfrac{2n}{ n-2}, \quad 2\leq \overline{r}_2< \infty,  \quad  &\mbox{ for }n\geq 2 \\[6pt]

2\leq \overline{q}_2\leq \infty, \quad \quad \ \,  2<\overline{r}_2<\infty, \quad   &\mbox{ for }n=1.
\end{cases}
\end{equation*}
 As before, we can find $0<\varkappa_2<1$ such that $\displaystyle{1 \over \overline{r}_2} + {n \over 2\overline{q}_2 }={n\over 4}(1+\varkappa_2)$ and 
  $ \widehat{q}_2= \overline{q}_2 (1+\varkappa_2)$ and $ \widehat{r}_2= \overline{r}_2 (1+\varkappa_2)$   verify   \eqref{intq,r} and by 
  \eqref{1.7LSU} 
 \begin{align*}
 \|u^{\alpha,k}\|^2_{\overline{q}_2, \overline{r}_2;Q_{t_1}^{\alpha,k}}&\leq \|u^{\alpha,k}\|^2_ {\widehat{q}_2, \widehat{r}_2;Q_{t_1}^{\alpha,k}} \, \|\chi_{Q_{t_1}^{\alpha,k}}\|^2_{{\widehat{q}_2\over \varkappa_2},{\widehat{r}_2\over \varkappa_2};Q^{\alpha,k}_{t_1}}\\ 
& \leq \|u^{\alpha,k}\|^2_ {\widehat{q}_2, \widehat{r}_2;Q_{t_1}^{\alpha,k}}\, \left(t_1^{{1\over \widehat{r}_2}}  |\Omega|^{{1\over\widehat{q}_2}} \right)^{2\varkappa_2}.
 \end{align*}
Hence, using \eqref{a1}, we get
\begin{equation} \label{2.8*}
\begin{split}
I_2  \leq 2\|\psi\|_{q_2,r_2;Q_{T}}&\, \bigg\{  \beta^2\|u^{\alpha,k}\|^2_{V_2^{1,0}(Q^{\alpha,k}_{t_1})} \big(t_1^{\frac{1}{\widehat{r}_2}}|\Omega|^{\frac{1}{\widehat{q}_2}}\big)^{2\varkappa_2}\\
&  +k^2\left[\int_0^{t_1}|A^{\alpha,k}(t)|^\frac{\widehat{r}_2}{\widehat{q}_2}\,dt\right]^{2(\varkappa_2 +1) \over \widehat{r}_2} \bigg\}.
\end{split}
\end{equation}
In case of $q_2=1$ and hence $\overline{q}_2=\infty,$  $n=1,$ we obtain 
\begin{equation} \label{2.8**}
\begin{split}
I_2  \leq 2\|\psi\|_{q_2,r_2;Q_{T}}&\, \left(  \beta^2\|u^{\alpha,k}\|^2_{V_2^{1,0}(Q^{\alpha,k}_{t_1})} \, t_1^{\frac{{2\varkappa_2}}{\widehat{r}_2}} +k^2\, t_1^{2(\varkappa_2 +1) \over \widehat{r}_2} \right).
\end{split}
\end{equation}  

Combine  \eqref{2.7*},  \eqref{2.8*} and  \eqref{2.5**}, we obtain 
\begin{equation}\label{2.9*}
\begin{split}
\|u^{\alpha,k}&\|^2_{V_2^{1,0}(Q_{t_1}^{\alpha,k})} \\ 
&\leq C\Biggl\{ \|u^{\alpha,k}\|_{V_2^{1,0}(Q_{t_1}^{\alpha,k})}^2 \Bigg[ \beta^2 \|u^\alpha\|^{\delta-2}_{q_1,r_1,Q_{t_1}^{\alpha,k}}  \big(t_1^{\frac{1}{\widehat{r}_1}} |\Omega |^{\frac{1}{\widehat{q}_1}} \big)^{2\varkappa_1} \\
&+ \beta^2 \|\psi\|_{q_2,r_2,Q_{T}}  \big(t_1^{\frac{1}{\widehat{r}_2}}|\Omega|^{\frac{1}{\widehat{q}_2}}\big)^{2\varkappa_2}\Bigg] \\
&+  \|u^\alpha\|^{\delta-2}_{q_1,r_1,Q_{t_1}^{\alpha,k}} \,k^2\left(\int_0^{t_1}\,|A^{\alpha,k}(t)|^{\frac{\widehat{r}_1}{\widehat{q}_1}}\,dt\right)^{\frac{2(\varkappa_1+1)}{\widehat{r}_1}}\\
& + k^2 \|\psi\|_{q_2,r_2,Q_{T}}\left(\int_0^{t_1}|A^{\alpha,k}(t)|^\frac{\widehat{r}_2}{\widehat{q}_2}\,dt\right)^{2(\varkappa_2 +1) \over \widehat{r}_2} \Biggl\}.
\end{split}
\end{equation}

Taking $t_1$ such small that 
$$
C \beta^2 \left(\|u^\alpha\|^{\delta-2}_{q_1,r_1,Q_{t_1}^{\alpha,k}}  \left(t_1^{\frac{1}{\widehat{r}_1}} |\Omega |^{\frac{1}{\widehat{q}_1}} \right)^{2\varkappa_1} 
+  \|\psi\|_{q_2,r_2,Q_{T}}  \left(t_1^{\frac{1}{\widehat{r}_2}}|\Omega|^{\frac{1}{\widehat{q}_2}}\right)^{2\varkappa_2}\right) < {1\over 2} 
$$
we obtain 
\begin{equation}\label{2.10}
\begin{split}
\|u^{\alpha,k}\|^2&_{V_{2}^{1,0}(Q_{t_1}^{\alpha,k})} \\
&\leq C\, k^2\Bigg[  \|u^\alpha\|^{\delta-2}_{q_1,r_1,Q_{t_1}^{\alpha,k}} \,\left(\int_0^{t_1}\,|A^{\alpha,k}(t)|^{\frac{\widehat{r}_1}{\widehat{q}_1}}\,dt\right)^{\frac{2(\varkappa_1+1)}{\widehat{r}_1}}\\
& +  \|\psi\|_{q_2,r_2,Q_{T}}\left(\int_0^{t_1}|A^{\alpha,k}(t)|^\frac{\widehat{r}_2}{\widehat{q}_2}\,dt\right)^{2(\varkappa_2 +1) \over \widehat{r}_2}\Bigg]. 
\end{split}
\end{equation}
Analogously in case of $n=1$.

Since  \eqref{2.10} is equivalent to \eqref{t1} in  Theorem \ref{T} we deduce that $u^\alpha$ is essentially bounded from above in a cylinder with  high $t_1$ which is  small enough. 
In order to extend this  result in the whole cylinder   we consider $u^\alpha$ successively in $Q_2=\Omega \times (t_1, 2t_1)$, $Q_3=\Omega \times (2t_1, 3t_1)$ and so on, covering in such way the whole cylinder $Q_T.$

Moreover, it is possible to estimate $u^\alpha(x,t)$ from below. To this aim, we apply the result just obtained to the function $\tilde{u}^\alpha(x,t)=-u^\alpha(x,t)$ which verifies equations similar to $u^\alpha(x,t)$, i.e.
\begin{equation*}
{\tilde{u}^\alpha}_t- \div \tilde{A}^\alpha(x,t,\tilde{\bu},D\tilde{\bu})+\tilde{b}^\alpha(x,t,\tilde{\bu},D\tilde{\bu})=0 \qquad (x,t)\in Q_T,
\end{equation*}
where 
$$\tilde{A}^\alpha(x,t,\tilde{\bu},D\tilde{\bu}) = - {A}^\alpha(x,t,-\tilde{\bu},-D\tilde{\bu})$$ 
and 
$$\tilde{b}^\alpha(x,t,\tilde{\bu},D\tilde{\bu})= -{b}^\alpha(x,t,-\tilde{\bu},-D\tilde{\bu}) $$ 
satisfy the conditions \eqref{1.8}-\eqref{cond2}. Thus, by  Theorem \ref{T}, we obtain 
$$\|u^\alpha \|_{\infty;Q_{T}}\leq K, $$ 
where  $K=K(M_0, k_0, q_i, r_i, n, \varkappa_i, | \Omega|,T, \|u^\alpha\|_{q_1,r_1;Q_T} ).$
\endproof

\subsection*{Acknowledgments.} 
\begin{small}
All the authors are members of \textit{INDAM - GNAMPA}.  
The research of L. Softova is partially supported by the project \textit{GNAMPA 2020 "Elliptic operators with unbounded and singular coefficients on weighted $L^p$ spaces".} \\
The research of E.A. Alfano is partially supported by the project \textit{ALPHA-MENTE,}  Lotto1/Ambito AV01.

\bigskip
 The authors are very indebted to the referee for the valuable remarks that improve  the exposition of the paper.
\end{small}

\end{document}